\documentclass[11pt]{article}
\usepackage{amsmath}
\usepackage{amssymb}
\usepackage{amsthm}
\usepackage[usenames]{color}
\usepackage{amscd}
\usepackage{indentfirst}

\def\1{{\bf 1}}

\def\N{{\Bbb N}}
\def\Z{{\Bbb Z}}

\def\C{{\Bbb C}}

\def\RE{\operatorname{Re}}

\def\lcm{\operatorname{lcm}}

\newtheorem{theorem}{Theorem}[section]
\newtheorem{corollary}{Corollary}[section]

\newtheorem{remark}{Remark}[section]
\newtheorem{lemma}{Lemma}[section]

\numberwithin{equation}{section}

\begin{document}

\title{\bf On the subgroups of finite Abelian groups of rank three}
\author{Mario Hampejs and L\'aszl\'o T\'oth\thanks{The second named author gratefully acknowledges support from the
Austrian Science Fund (FWF) under the project Nr. M1376-N18.}
}
\date{}
\maketitle

\centerline{{\sl Annales Univ. Sci. Budapest., Sect. Comp.} {\bf 39}
(2013), 111--124}

\begin{abstract} We describe the subgroups of the group $\Z_m \times \Z_n
\times \Z_r$ and derive a simple formula for the total number
$s(m,n,r)$ of the subgroups, where $m,n,r$ are arbitrary positive
integers. An asymptotic formula for the function $n\mapsto s(n,n,n)$
is also deduced.
\end{abstract}




\section{Introduction}

Throughout the paper we use the notation: $\N=\{1,2,\ldots\}$,
$\N_0=\{0,1,2,\ldots\}$, $\Z_m$ is the additive group of residue
classes modulo $m$, $\phi$ is Euler's totient function, $\tau(n)$ is
the number of divisors of $n$, $\zeta$ is the Riemann zeta function.

For an arbitrary finite Abelian group $G$ of order $\# G$ let $s(G)$
denote the total number of its subgroups. It is known that the
problem of counting the subgroups of $G$ reduces to $p$-groups. More
precisely, let $\# G= p_1^{a_1}\cdots p_r^{a_r}$ be the prime power
factorization of $\# G$ and let $G=G_1\times \cdots \times G_r$ be
the primary decomposition of $G$, where $\# G_i=p_i^{a_i}$ ($1\le
i\le r$). Then
\begin{equation*}
s(G)=s(G_1)\cdots s(G_r),
\end{equation*}
which follows from the properties of the subgroup lattice of $G$.
See, e.g., R.~Schmidt \cite{Sch1994} and M.~Suzuki \cite{Suz1951}.

Now let $G_{(p)}$ be a $p$-group of type $\lambda=
(\lambda_1,\ldots,\lambda_r)$, with $\lambda_1\ge \ldots \ge
\lambda_r\ge 1$, where $\lambda$ is a partition of
$|\lambda|=\lambda_1+\ldots +\lambda_r$. Formulas for the number
$s_{\mu}(G_{(p)})$ of subgroups of type $\mu$ ($\mu \subseteq
\lambda$) of $G_{(p)}$ were established by several authors, see
G.~Birkhoff \cite{Bir1934}, S.~Delsarte \cite{Del1948},
P.~E.~Dyubyuk \cite{Dyu1948}, Y.~Yeh \cite{Yeh1948}. One of these
formulas is given, in terms of the Gaussian coefficients $\left[{r
\atop k} \right]_p= \prod_{i=1}^k \frac{ p^{r-k+i}-1}{p^i-1}$ by
\begin{equation} \label{subgr_type}
s_{\mu}(G_{(p)})= \prod_{j=1}^{\lambda_1}
p^{\mu'_{j+1}(\lambda'_j-\mu'_j)} \left[{\lambda'_j-\mu'_{j+1} \atop
\mu'_j-\mu'_{j+1}} \right]_p,
\end{equation}
where $\lambda'$ and $\mu'$ are the conjugates (according to the
Ferrers diagrams) of $\lambda$ and $\mu$, respectively. Hence
$s_{\mu}(G_{(p)})$ is a polynomial in $p$, with integer
coefficients, depending only on $\lambda$ and $\mu$ (it is a sum of
Hall polynomials). Therefore, the number of the subgroups of order
$p^k$ ($0\le k\le |\lambda|$) of $G_{(p)}$ is
\begin{equation*}
s_{p^k}(G_{(p)}) = \sum_{\substack{\mu\subseteq \lambda\\ |\mu|=k}}
s_{\mu}(G_{(p)})
\end{equation*}
and the total number of subgroups is given by $s(G_{(p)}) =
\sum_{0\le k\le |\lambda|} s_{\mu}(G_{(p)})$. See the monograph of
M.~L.~Butler \cite{But1994} for a detailed discussion of formula
\eqref{subgr_type} and of related results of which proofs are
combinatorial and linear algebraic in nature.

Another general formula for the total number of subgroups of a
$p$-group of arbitrary rank, obtained by combinatorial arguments
using divisor functions of matrices  was given by G.~Bhowmik
\cite{Bho1996}. However, it is rather complicate to apply the above
formulas or that of \cite{Bho1996} to compute numerically the total
number of subgroups (of a given order) of a $p$-group. Also it is
difficult to find the coefficients of the polynomials in $p$
representing the number of subgroups (of a given order) of a
$p$-group, even in the case of rank two or three.

There are other tools which can be used to derive explicit formulas
for the total number of subgroups in the case of $p$-groups of rank
two. Namely, Goursat's lemma for groups was applied by
G.~C\u{a}lug\u{a}reanu \cite{Cal2004} and J.~Petrillo
\cite{Pet2011}, and the concept of the fundamental group lattice was
used by M.~T\u{a}rn\u{a}uceanu \cite{Tar2007,Tar2010}. In the paper
\cite{HHTW2012} the subgroups of $\Z_m \times \Z_n$ were
investigated, where $m,n\in \N$ are arbitrary, and the following
compact formula was deduced. The total number $s(m,n)$ of subgroups
of $\Z_m \times \Z_n$ is given by
\begin{equation} \label{form_2}
s(m,n)=\sum_{a\mid m, b\mid n} \gcd(a,b).
\end{equation}

Consider now the case of $p$-groups of rank three. It is well known
that for every prime $p$ the elementary Abelian group $(\Z_p)^3$ can
be considered as a three dimensional linear space over the Galois
field $GF(p)$. Its $k$-dimensional subspaces are exactly the
subgroups of order $p^k$  and the number of these subgroups is given
by the Gaussian coefficients $\left[{3 \atop k} \right]_p$ ($0\le
k\le 3$). The total number of subgroups of $(\Z_p)^3$ is $s(p)=
\sum_{k=0}^3 \left[{3 \atop k} \right]_p= 2(p^2+p+2)$. Similar
considerations hold also for the elementary Abelian groups
$(\Z_p)^r$ with $r\in \N$, cf. \cite{Aig2007,But1994,
Cal2004,Tar2010}.

It seems that, excepting the case of $(\Z_p)^3$ no simple general
formulas are known in the literature to generate the subgroups and
to compute the number of the subgroups of an Abelian group of rank
three. We refer here also to the paper of R.~Remak \cite[Sect.\
2]{Rem1931}, concerning a more general case, namely the direct
product of three finite groups, but where some $156$ equations are
given to describe the subgroups.

In this paper we investigate the subgroups of $p$-groups of rank
three. In fact, we consider the group $\Gamma:=\Z_m \times \Z_n
\times \Z_r$, where $m,n,r$ are arbitrary positive integers,
describe its subgroups and derive a simple formula for the total
number $s(m,n,r)$ of subgroups of $\Gamma$. We also deduce an
asymptotic formula for the function $n\mapsto s(n):=s(n,n,n)$.

Our approach is elementary, different from those quoted above, using
only simple group-theoretic and number-theoretic arguments. The main
results are given in Section \ref{Sect_results}, while their proofs
are presented in Section \ref{Sect_proofs}. Section \ref{Sect_val}
includes tables with numerical values and formulae regarding
$s(m,n,r)$.

We also remark that the number $c(m,n,r)$ of cyclic subgroups of
$\Gamma=\Z_m \times \Z_n \times \Z_r$ is given by
\begin{equation*} \label{total_number_cyclic_subgroups_var_1}[section]
c(m,n,r) = \sum_{a\mid m, b\mid n, c\mid r} \frac{\phi(a)
\phi(b)\phi(c)}{\phi(\lcm(a,b,c))},
\end{equation*}
see \cite{Tot2011,Tot2012}. The functions $(m,n,r)\mapsto s(m,n,r)$
and $(m,n,r)\mapsto c(m,n,r)$ are multiplicative functions of three
variables (cf. e.g., \cite{Tot2011} for this notion). The function
$n\mapsto s(n)$ is a multiplicative function of a single variable,
it is the sequence \cite[item A064803]{OEIS}.

\section{Results} \label{Sect_results}

Our first result is concerning the representation of the subgroups
of $\Gamma$.

\begin{theorem} \label{Th_repr} Let $m,n,r\in \N$. The subgroups of
the group $\Gamma=\Z_m \times \Z_n \times \Z_r$ can be represented
as follows.

(i) Choose $a,b,c\in \N$ such that $a\mid m, b\mid n, c\mid r$.

(ii) Compute $A:= \gcd(a,n/b)$, $B:=\gcd(b,r/c)$, $C:=\gcd(a,r/c)$.

(iii) Compute
\begin{equation*}
X:= \frac{ABC}{\gcd(a(r/c),ABC)}.
\end{equation*}

(iv) Let $s:=at/A$, where $0\le t\le A-1$.

(v) Let
\begin{equation*}
v:=\frac{bX}{B\gcd(t,X)}w, \text{ where } 0\le w\le B\gcd(t,X)/X-1.
\end{equation*}

(vi) Find a solution $u_0$ of the linear congruence
\begin{equation*}
(r/c)u \equiv rvs/(bc) \quad  \text{(mod $a$)}.
\end{equation*}

(vii) Let $u:=u_0+az/C$, where $0\le z\le C-1$.

(viii) Consider
\begin{equation*}
U_{a,b,c,t,w,z}:= \langle (a,0,0),(s,b,0),(u,v,c) \rangle
\end{equation*}
\begin{equation*}
= \{(ia+js+ku, jb+kv,kc): 0\le i\le n/a-1, 0\le j\le n/b-1, 0\le
k\le n/c-1 \}.
\end{equation*}

Then $U_{a,b,c,t,w,z}$ is a subgroup of order $mnr/(abc)$ of
$\Gamma$. Moreover, there is a bijection between the set of
sextuples $(a,b,c,t,w,z)$ satisfying the conditions (i)-(viii) and
the set of subgroups of $\Gamma$.
\end{theorem}

Let $P(n):=\sum_{k=1}^n \gcd(k,n)=\sum_{d\mid n} d\phi(n/d)$ be the
gcd-sum function. Note that the function $P$ is multiplicative and
\begin{equation} \label{form_P}
P(p^\nu)= (\nu+1)p^\nu - \nu p^{\nu -1}
\end{equation}
for every prime power $p^\nu$ ($\nu \in \N$). See \cite{Tot2010}.

Next we give a formula for the number of subgroups of $\Gamma$.

\begin{theorem} \label{Th_number}  For every $m,n,r\in \N$ the total number of subgroups
of the group $\Z_m \times \Z_n \times \Z_r$ is given by
\begin{equation} \label{form_3}
s(m,n,r)= \sum_{a\mid m, b\mid n, c\mid r} \frac{ABC}{X^2} P(X),
\end{equation}
with the notation of Theorem \ref{Th_repr}.

The number of subgroups of order $\delta$ ($\delta \mid mnr$) is
given by \eqref{form_3} with the additional condition that the
summation is subject to $abc= mnr/\delta$.
\end{theorem}

If one of $m,n,r$ is $1$, then formula \eqref{form_3} reduces to
\eqref{form_2}.

\begin{corollary} \label{Cor} For every prime $p$ and every $\nu_1,\nu_2,\nu_3 \in \N$,
$s(p^{\nu_1},p^{\nu_2},p^{\nu_3})$ is a polynomial in $p$ with
integer coefficients.

In particular, for every $\nu \in \N$, $s(p^{\nu})$ is a polynomial
in $p$ of degree $2\nu$, having the leading coefficient $\nu+1$.
\end{corollary}

See Section \ref{Sect_val}, Table 2 for the polynomials $s(p^{\nu})$
with $1\le \nu \le 10$.

\begin{remark} Actually, for every $\nu \in \N$,
\begin{equation} \label{general_form}
s(p^{\nu})= \sum_{j=0}^{2\nu} \left(\nu -
\left[\frac{j-1}{2}\right]\right) \left(2j-
\left[\frac{j-1}{2}\right]\right)p^{2\nu-j},
\end{equation}
where $[x]$ denotes the integer part of $x$. A proof of
\eqref{general_form} will be presented elsewhere.
\end{remark}

The asymptotic behavior of the function $n\mapsto s(n)$ is related
to Dirichlet's divisor problem. Let $\theta$ be the number such that
\begin{equation} \label{Dirichlet_divisor}
\sum_{n\le x} \tau(n) = x\log x + (2\gamma-1)x + {\cal
O}(x^{\theta+\varepsilon}),
\end{equation}
for every $\varepsilon >0$, where $\gamma$ is the Euler-Mascheroni
constant. It is known that $1/4\le \theta \le 131/416 \approx
0.3149$, where the upper bound, the best up to date, is the result
of M.~N.~Huxley \cite{Hux2003}. The following asymptotic formula
holds. Define the multiplicative function $h$ by
\begin{equation} \label{def_h}
s(n)= \sum_{d\mid n} d^2\tau(d) h(n/d) \quad (n\in \N)
\end{equation}
and let $H(z)= \sum_{n=1}^{\infty} h(n)n^{-z}$ be the Dirichlet
series of $h$.

\begin{theorem} \label{Th_asymp} For every $\varepsilon >0$,
\begin{equation} \label{asymp_s}
\sum_{n\le x} s(n) = \frac{x^3}{3} \left(H(3)( \log x+ 2\gamma-1)+
H'(3) \right) + {\cal O}(x^{2+\theta+ \varepsilon}),
\end{equation}
where $H'$ is the derivative of $H$.
\end{theorem}

\begin{remark} It follows from \eqref{general_form} and \eqref{def_h} by induction
that $h(p^{\nu})= (3\nu-1)p+3\nu+1$ for
every prime power $p^{\nu}$ ($\nu \in \N$). Hence,
\begin{equation} \label{def_series_H}
H(z)=\zeta^2(z)\prod_p \left(1+\frac{2}{p^{z-1}} + \frac{2}{p^z}+
\frac{1}{p^{2z-1}}\right)
\end{equation}
for $\RE(z)>2$. In particular,
\begin{equation*}
H(3)=\zeta^2(3)\prod_p \left(1+\frac{2}{p^2} + \frac{2}{p^3}+
\frac{1}{p^5}\right).
\end{equation*}
\end{remark}

\section{Proofs} \label{Sect_proofs}

We need the next general result regarding the subgroups of the group
$G\times \Z_q$, where $(G,+)$ is an arbitrary finite Abelian group.
For a subgroup $H$ of $G$ (notation $H\le G$) consider the
congruence relation $\varrho_H$ on $G$ defined for $x,x'\in G$ by $x
\varrho_H \, x'$ if $x-x'\in H$.

\begin{lemma} \label{Lemma_1} For a finite Abelian group $(G,+)$ and
$q\in \N$ let
\begin{equation*} I_{G,q}:= \{(H,\alpha,d): H\le G,
\alpha \in {\cal S}_H, d\mid q \ {\text and } \ (q/d)\alpha \in H\},
\end{equation*}
where ${\cal S}_H$ is a complete system of representants of the
equivalence classes determined by $\varrho_H$. For $(H,\alpha,d)\in
I_{G,q}$ define
\begin{equation*}
V_{H,\alpha,d}:= \{(k\alpha+\beta,kd): 0\le k\le q/d-1, \beta \in
H\}.
\end{equation*}

Then $V_{H,\alpha,d}$ is a subgroup of order $(q/d)\# H$ of $G\times
\Z_q$ and the map $(H,\alpha,d) \mapsto V_{H,\alpha,d}$ is a
bijection between the set $I_{G,q}$ and the set of subgroups of
$G\times \Z_q$.
\end{lemma}

\proof Let $V$ be a subgroup of $G\times \Z_q$. Consider the natural
projection $\pi_2:G\times \Z_q \to \Z_q$ given by $\pi_2(x,y)=y$.
Then $\pi_2(V)$ is a subgroup of $\Z_q$ and there is a unique
divisor $d$ of $q$ such that $\pi_2(V)=\langle d \rangle:=\{kd: 0\le
k\le q/d-1\}$. Let $\alpha \in G$ such that $(\alpha,d)\in V$.

Furthermore, consider the natural inclusion $\iota_1:G\to G\times
\Z_q$ given by $\iota_1(x)=(x,0)$. Then $\iota_1^{-1}(V)=H$ is a
subgroup of $G$. We show that $V=\{(k\alpha+\beta,kd): k\in \Z,
\beta \in H\}$. Indeed, for every $k\in \Z$ and $\beta \in H$,
$(k\alpha+\beta,kd)= k(\alpha,d)+(\beta,0)\in V$. On the other hand,
for every $(u,v)\in V$ one has $v\in \pi_2(V)$ and hence there is
$k\in \Z$ such that $v=kd$. We obtain $(u-k\alpha,0)= (u,v)-
k(\alpha,d)\in V$, thus $\beta:=u-k\alpha \in \iota_1^{-1}(V)=H$.

Here a necessary condition is that $(q/d)\alpha \in H$ (obtained for
$k=q/d$, $\beta=0$). Clearly, if this is verified, then for the
above representation of $V$ it is enough to take the values $0\le
k\le q/d-1$.

Conversely, every $(H,\alpha,d)\in I_{G,q}$ generates a subgroup
$V_{H,\alpha,d}$ of order $(q/d)\# H$ of $G\times \Z_q$.
Furthermore, for fixed $H\le G$ and $d\mid q$ we have
$V_{H,\alpha,d} = H_{H,\alpha',d}$ if and only if $\alpha \varrho_H
\, \alpha'$. This completes the proof. \qed

In the case $G=\Z_m$ (and with $q=n$) Lemma \ref{Lemma_1} was given
in \cite[Th.\ 1]{HHTW2012} and it can be stated as follows:

\begin{lemma} \label{Lemma_2} For every $m,n\in \N$ let
\begin{equation*}
I_{m,n}:=\{(a,b,s)\in \N^2\times \N_0: a\mid m, b\mid n, 0\le s\le
a-1 \ {\text and } \ a\mid (n/b)s \}
\end{equation*}
and for $(a,b,s)\in I_{m,n}$ define
\begin{equation} \label{B}
V_{a,b,s}:= \langle (a,0), (b,s) \rangle
\end{equation}
\begin{equation*}
=\{(ia+js,jb): 0\le i\le m/a-1, 0\le j\le n/b-1\}.
\end{equation*}

Then $V_{a,b,s}$ is a subgroup of order $\frac{mn}{ab}$ of
$\Z_m\times \Z_n$ and the map $(a,b,s)\mapsto V_{a,b,s}$ is a
bijection between the set $I_{m,n}$ and the set of subgroups of
$\Z_m \times \Z_n$.
\end{lemma}

Note that $a\mid (n/b)s$ holds if and only if $a/\gcd(a,n/b) \mid
s$. That is, for $s\in I_{m,n}$ we have
\begin{equation} \label{s}
s= \frac{at}{A}, \quad 0\le t\le A-1,
\end{equation}
where $A=\gcd(a,n/b)$, notation given in Theorem \ref{Th_repr}. This
leads quickly to formula \eqref{form_2} regarding the number
$s(m,n)$ of subgroups of $\Z_m \times \Z_n$, namely
\begin{equation*}
s(m,n)= \sum_{a\mid m, b\mid n} \ \sum_{0\le t\le A-1} 1
=\sum_{a\mid m, b\mid n} \gcd(a,b).
\end{equation*}

\proof (for Theorem \ref{Th_repr}) Apply Lemma \ref{Lemma_1} for
$G=\Z_m \times \Z_n$ and with $q=r$. For the subgroups $V=V_{a,b,s}$
given by Lemma \ref{Lemma_2} a complete system of representants of
the equivalence classes determined by $\varrho_V$ is ${\cal
S}_{a,b}=\{0,1,\ldots, a-1\} \times \{0,1,\ldots,b-1\}$. Indeed, the
elements of ${\cal S}_{a,b}$ are pairwise incongruent with respect
to $V$, and for every $(x,y)\in \Z_m \times \Z_n$ there is a unique
$(x',y')\in {\cal S}_{a,b}$ such that $(x,y)-(x',y')\in V$. Namely,
let
\begin{equation*}
(x_1,y_1)= (x,y)- \lfloor y/b \rfloor (s,b), \qquad (x',y')=
(x_1,y_1)- \lfloor x_1/a \rfloor (a,0).
\end{equation*}

We obtain that the subgroups of $\Z_m \times \Z_n \times \Z_r$ are
of the form
\begin{equation*}
U=U_{H,\alpha,c}= \{(k\alpha+\beta,kc): 0\le k\le r/c-1, \beta \in
V\},
\end{equation*}
where $c\mid r$ and $\alpha =(u,v)\in {\cal S}_{a,b}$ such that
$(r/c)\alpha \in V$.

Now using \eqref{B} we deduce
\begin{equation*}
U=U_{a,b,s,\alpha,c}
\end{equation*}
\begin{equation*}
= \{(ia+js+ku, jb+kv,kc): 0\le i\le n/a-1, 0\le j\le n/b-1, 0\le
k\le n/c-1 \},
\end{equation*}
where \eqref{s} holds and $(r/c)(u,v)\in V$. From the latter
condition we deduce that there are $i_0,j_0$ such that
\begin{equation} \label{cond}
(r/c)u =i_0a+j_0s, \qquad  (r/c)v=j_0b.
\end{equation}

The second condition of \eqref{cond} holds if $b\mid (r/c)v$, that
is $b/\gcd(b,r/c) \mid v$. Let
\begin{equation} \label{v}
v= \frac{bv_1}{B}, \quad 0\le v_1\le B-1,
\end{equation}
where $B=\gcd(b,r/c)$. Also, $j_0=rv/(bc)$ and inserting this into
the first equation of \eqref{cond} we obtain $(r/c)u \equiv
rvs/(bc)$ (mod $a$). This linear congruence in $u$ has a solution
$u_0$ if and only if
\begin{equation} \label{new_cond}
\gcd(a,r/c) \mid \frac{rvs}{bc}
\end{equation}
and its all solutions are $u=u_0+az/C$ with $0\le z\le C-1$ with
$C=\gcd(a,r/c)$.

Substituting \eqref{v} and \eqref{s} into \eqref{new_cond} we obtain
\begin{equation*}
\gcd(a,r/c) \mid \frac{rab}{\gcd(ab,n)\gcd(bc,r)} v_1t,
\end{equation*}
that is
\begin{equation*}
\gcd(ab,n) \gcd(ac,r) \gcd(bc,r) \mid abcr v_1t,
\end{equation*}
equivalent to
\begin{equation*}
\frac{\gcd(ab,n)\gcd(ac,r) \gcd(bc,r)}{\gcd(abcr,
\gcd(ab,n)\gcd(ac,r)\gcd(bc,r))} \mid v_1t,
\end{equation*}
and to
\begin{equation} \label{cond_X}
X\mid v_1t,
\end{equation}
where $X$ is defined in the statement of Theorem \ref{Th_repr}. Note
that $X\mid B$ (indeed, $A\mid a$, $C\mid (r/c)$ and the property
follows from $X=B/\gcd((a/A)(r/c)/C,B$).

Let $t$ be fixed. We obtain from \eqref{cond_X} that $v_1$ is of the
form $v_1=Xw/\gcd(t,X)$, where $0\le w \le B\gcd(t,X)/X-1$. Also,
from \eqref{v}, $v=bXw/B\gcd(t,X)$. Collecting the conditions on
$a,b,c,t,w,z$ in terms of $A,B,C,X$ finishes the proof. \qed

\proof (for Theorem \ref{Th_number}) According to Theorem
\ref{Th_repr} the number of subgroups of $\Gamma$ is
\begin{equation*}
s(m,n,r)= \sum_{a\mid m, b\mid n, c\mid r} \ \sum_{0\le t\le A-1} \
\sum_{0\le w\le B\gcd(t,X)/X-1} \ \sum_{0\le z\le C-1} 1
\end{equation*}
\begin{equation*}
= \sum_{a\mid m, b\mid n, c\mid r} C \sum_{0\le t\le A-1}
\frac{B}{X} \gcd(t,X) = \sum_{a\mid m, b\mid n, c\mid r}
\frac{BC}{X} \sum_{1\le t\le A} \gcd(t,X).
\end{equation*}

Here $X\mid A$ (similar to $X\mid B$ shown above), hence the inner
sum is $(A/X)P(X)$ and we obtain the formula \eqref{form_3}.

In the case of subgroups of order $\delta$ use that the order of
$U_{a,b,c,t,w,z}$ is $mnr/(abc)$, according to Theorem
\ref{Th_repr}. \qed

\proof (for Corollary \ref{Cor}) If $m=p^{\nu_1}$, $n=p^{\nu_2}$,
$r=p^{\nu_3}$, then for each term of the sum \eqref{form_3} all of
$A,B,C,X$ and $ABC/X^2$ are of form $p^{\nu}$ with some $\nu \in
\N_0$ ($X\mid A$ and $X\mid B$, cf. the proof of Theorem
\ref{Th_repr}). Using the formula \eqref{form_P} for $P(p^\nu)$ we
deduce that $s(p^{\nu_1}, p^{\nu_2}, p^{\nu_3})$ is a polynomial in
$p$ with integer coefficients.

In the case $\nu_1=\nu_2=\nu_3=\nu$, by writing explicitly the terms
of the sum \eqref{form_3} for various choices of $a,b,c\mid
p^{\nu}$, we deduce that the maximal exponent of $p$ is $2\nu$,
which is obtained exactly for $a=p^{\nu}$, $b=p^{\lambda}$ ($0\le
\lambda \le \nu$) and $c=1$. \qed

\proof (for Theorem \ref{Th_asymp}) It follows from
\eqref{def_series_H} that the abscissa of absolute convergence of
$H(z)$ is $2$. But \eqref{def_series_H} is a consequence of the
formula \eqref{general_form}, not proved in the present paper. For
this reason we show here by different arguments that $H(z)$ is
absolutely convergent for every $z\in \C$ with $\RE(z)>
9/4+\varepsilon$, which is sufficient to establish the asymptotic
formula.

According to \eqref{def_h}, the function $s$ can be expressed in
terms of the Dirichlet convolution $*$ as $s=E_2\tau*h$, where
$E_2(n)=n^2$ ($n\in \N$). Therefore, $h=s*(\mu*\mu)E_2$, $\mu$
denoting the M\"obius function. We obtain that $h(p)=2p+4$,
$h(p^2)=5p+7$, $h(p^3)=8p+10$ and
\begin{equation} \label{express_h}
h(p^{\nu})= s(p^{\nu})- 2p^2s(p^{\nu-1}) + p^4s(p^{\nu-2}) \quad (\nu \ge 2).
\end{equation}

For the gcd-sum function $P$ one has $P(n)\le n\tau(n)$ ($n\in \N$), cf. \eqref{form_P}. Hence
\begin{equation*}
s(n)=\sum_{a,b,c\mid n} \gcd(a(r/c),ABC) P(X)/X \le
\sum_{a,b,c\mid n} a(r/c) \tau(X)
\end{equation*}
\begin{equation*}
\le n^2 \tau(n) \sum_{a,b,c\mid n} 1 = n^2 (\tau(n))^4
\end{equation*}
for every $n\in \N$ and for every prime power $p^{\nu}$ ($\nu \in \N$),
\begin{equation} \label{ineq_s}
s(p^{\nu}) \le p^{2\nu} (\nu+1)^4.
\end{equation}

Now from \eqref{express_h} and \eqref{ineq_s} we deduce that for every prime power $p^{\nu}$ ($\nu \ge 2$),
\begin{equation} \label{ineq_h}
0< h(p^{\nu}) \le 2 p^{2\nu} (\nu+1)^4.
\end{equation}

From the Euler product formula
\begin{equation*}
H(z) = \prod_p \left(1+\frac{2p+4}{p^z}+ \frac{5p+7}{p^{2z}}+ \frac{8p+10}{p^{3z}}+ \sum_{\nu=4}^{\infty}
\frac{h(p^{\nu})}{p^{\nu z}} \right)
\end{equation*}
and from \eqref{ineq_h} we obtain that $H(z)$ is absolutely
convergent for $z\in \C$ with $4(\RE(z)-2-\varepsilon)
>1$, i.e., for $\RE (z) >9/4+\varepsilon$ with an arbitrary $\varepsilon >0$.

Furthermore, by partial summation we obtain from \eqref{Dirichlet_divisor} that
\begin{equation} \label{n2_tau}
\sum_{n\le x} n^2\tau(n) = \frac1{3} x^3 \log x + \frac1{3}
\left(2\gamma-\frac1{3}\right)x^3 + {\cal
O}(x^{2+\theta+\varepsilon}).
\end{equation}

Now
\begin{equation*}
\sum_{n\le x} s(n)= \sum_{d\le x} h(d) \sum_{e\le x/d} e^2\tau(e),
\end{equation*}
and inserting \eqref{n2_tau} we get
\begin{equation*}
\sum_{n\le x} s(n) = \frac{x^3\log x}{3} \sum_{d\le x}
\frac{h(d)}{d^3} - \frac{x^3}{3} \sum_{d\le x} \frac{h(d)\log
d}{d^3} +\frac{x^3}{3} \left(2\gamma-\frac1{3}\right) \sum_{d\le x}
\frac{h(d)}{d^3}
\end{equation*}
\begin{equation*}
+ {\cal O}\left(x^{2+\theta+\varepsilon} \sum_{d\le x} \frac{|h(d)|}{d^{2+\theta+\varepsilon}}\right),
\end{equation*}
where the last term is ${\cal O}\left(x^{2+\theta+
\varepsilon}\right)$. This gives the asymptotic formula
\eqref{asymp_s}. \qed

\section{Tables} \label{Sect_val}

The computations were performed using the software Mathematica.

\medskip \medskip

\centerline{Table 1. Values of $s(n)$ for $1\le n\le 50$}
\vskip-4mm
\[
\vbox{\offinterlineskip \hrule \halign{ \strut
\vrule \hfill $\ # \ $ \hfill
& \vrule \hfill $\ # \ $ \hfill
& \vrule \vrule  \hfill $\ # \ $ \hfill
& \vrule \hfill $\ # \ $ \hfill
& \vrule \vrule \hfill $\ # \ $  \hfill
& \vrule \hfill $\ # \ $ \hfill
& \vrule \vrule \hfill $\ # \ $ \hfill
& \vrule \hfill $\ # \ $ \hfill
& \vrule \vrule \hfill $\ # \ $ \hfill
& \vrule \hfill $\ # \ $  \hfill \vrule \cr
 n & s(n) & n & s(n) & n & s(n) & n & s(n) & n & s(n) \ \cr \noalign{\hrule}
 1 & 1 & 11 & 268 & 21  & 3248  & 31  & 1988  & 41 & 3448 \cr \noalign{\hrule}
 2 & 16 & 12 & 3612 & 22  & 4288  & 32  & 22308  & 42 & 51968 \cr \noalign{\hrule}
 3 & 28 & 13 & 368 & 23  & 1108  & 33  & 7504  & 43 & 3788 \cr \noalign{\hrule}
 4 & 129 & 14 & 1856 & 24  & 22456  & 34  & 9856  & 44 & 34572 \cr \noalign{\hrule}
 5 & 64 & 15 & 1792 &  25 & 2607  &  35 & 7424  & 45 & 28480 \cr \noalign{\hrule}
 6 & 448 & 16 & 4387 &  26 & 5888  &  36 & 57405  & 46 & 17728 \cr \noalign{\hrule}
 7 & 116 & 17 & 616 &  27 & 5776  &  37 & 2816  & 47 & 4516 \cr \noalign{\hrule}
 8 & 802 & 18 & 7120 &  28 & 14964  &  38 & 12224  & 48 & 122836 \cr \noalign{\hrule}
 9 & 445 & 19 & 764 &  29 & 1744  &  39 & 10304  &  49 & 9009 \cr \noalign{\hrule}
 10 &  1024 &  20 & 8256 &  30 & 28672  &  40 & 51328  &  50 & 41712  \cr} \hrule}
\]

\medskip
\newpage

\centerline{Table 2. Values of $s(p^{\nu})$ for $1\le \nu \le 10$}
\vskip-4mm
\[
\vbox{\offinterlineskip  \hrule \halign{ \strut
\vrule \hfill $\ # \ $
& \vrule  $\ #  $ \hfill  \vrule  \cr
\nu & \hfill s(p^{\nu}) \ \cr \noalign{\hrule}
 1 & 4+2 p+2 p^2 \cr \noalign{\hrule}
 2 & 7+5 p+8 p^2+4 p^3+3 p^4 \cr \noalign{\hrule}
 3 & 10+8 p+14 p^2+10 p^3+12 p^4+6 p^5+4 p^6 \cr \noalign{\hrule}
 4 & 13+11 p+20 p^2+16 p^3+21 p^4+15 p^5+16 p^6+8 p^7+5 p^8 \cr \noalign{\hrule}
 5 & 16+14 p+26 p^2+22 p^3+30 p^4+24 p^5+28 p^6+20 p^7+20 p^8 \cr
   & +10 p^9+6 p^{10} \cr \noalign{\hrule}
 6 & 19+17 p+32 p^2+28 p^3+39 p^4+33 p^5+40 p^6+32 p^7+35 p^8 \cr
   & +25 p^9+24 p^{10}+12 p^{11}+7 p^{12} \cr \noalign{\hrule}
 7 & 22+20 p+38 p^2+34 p^3+48 p^4+42 p^5+52 p^6+44 p^7+50 p^8 \cr
   & + 40 p^9+42 p^{10}+30 p^{11}+28 p^{12} +14 p^{13}+8 p^{14} \cr \noalign{\hrule}
 8 & 25+ 23 p +44 p^2 + 40 p^3 +57 p^4 + 51 p^5 + 64 p^6 + 56 p^7 +65 p^8 \cr
   &  +55p^9 + 60 p^{10} +48 p^{11} +49 p^{12} + 35 p^{13} + 32 p^{14} +16p^{15}+ 9p^{16} \cr \noalign{\hrule}
 9 & 28+ 26 p + 50p^2 + 46 p^3+66 p^4+ 60p^5 + 76p^6 + 68p^7 + 80p^8 \cr
   & + 70p^9 + 78 p^{10} + 66 p^{11}+ 70 p^{12} +56 p^13 + 56 p^{14} + 40 p^{15} + 36 p^{16} \cr
   & + 18 p^{17} +10 p^{18} \cr \noalign{\hrule}
10 & 31 + 29p + 56p^2 + 52p^3 + 75 p^4 + 69 p^5 + 88 p^6 + 80 p^7+
95 p^8 \cr
   & + 85 p^9+ 96 p^{10} + 84 p^{11} + 91 p^{12} + 77 p^{13} + 80 p^{14} + 64p^{15} + 63p^{16}\cr
   & + 45p^{17} +40p^{18} +20 p^{19} + 11p^{20}  \cr} \hrule}
\]

\medskip
\newpage

\centerline{Table 3. Values of $s(p^{\nu_1},p^{\nu_2},p^{\nu_3})$
for $1\le \nu_1 \le \nu_2 \le \nu_3\le 4$} \vskip-4mm
\[
\vbox{\offinterlineskip  \hrule \halign{ \strut \vrule \hfill $\ # \
$ & \vrule \hfill $\ # \ $ \hfill & \vrule  \hfill $\ # \ $ \hfill &
\vrule $\ # \ $ \hfill \vrule  \cr \nu_1 & \nu_2 & \nu_3 & \hfill
s(p^{\nu_1},p^{\nu_2},p^{\nu_3})\ \cr \noalign{\hrule}
 1 & 1 & 1 & 4+2 p+2 p^2 \cr \noalign{\hrule}
 1 & 1 & 2 & 5+3 p+4 p^2 \cr \noalign{\hrule}
 1 & 2 & 2 & 6+4 p+6 p^2+2 p^3 \cr \noalign{\hrule}
 2 & 2 & 2 & 7+5 p+8 p^2+4 p^3+3 p^4 \cr \noalign{\hrule}
 1 & 1 & 3 & 6+4 p+6 p^2 \cr \noalign{\hrule}
 1 & 2 & 3 & 7+5 p+8 p^2+4 p^3 \cr \noalign{\hrule}
 2 & 2 & 3 & 8+6 p+10 p^2+6 p^3+6 p^4 \cr \noalign{\hrule}
 1 & 3 & 3 & 8+6 p+10 p^2+6 p^3+2 p^4 \cr \noalign{\hrule}
 2 & 3 & 3 & 9+7 p+12 p^2+8 p^3+9 p^4+3 p^5 \cr \noalign{\hrule}
 3 & 3 & 3 & 10+8 p+14 p^2+10 p^3+12 p^4+6 p^5+4 p^6 \cr \noalign{\hrule}
 1 & 1 & 4 & 7+5 p+8 p^2 \cr \noalign{\hrule}
 1 & 2 & 4 & 8+6 p+10 p^2+6 p^3 \cr \noalign{\hrule}
 2 & 2 & 4 & 9+7 p+12 p^2+8 p^3+9 p^4 \cr \noalign{\hrule}
 1 & 3 & 4 & 9+7 p+12 p^2+8 p^3+4 p^4 \cr \noalign{\hrule}
 2 & 3 & 4 & 10+8 p+14 p^2+10 p^3+12 p^4+6 p^5 \cr \noalign{\hrule}
 3 & 3 & 4 & 11+9 p+16 p^2+12 p^3+15 p^4+9 p^5+8 p^6 \cr \noalign{\hrule}
 1 & 4 & 4 & 10+8 p+14 p^2+10 p^3+6 p^4+2 p^5 \cr \noalign{\hrule}
 2 & 4 & 4 & 11+9 p+16 p^2+12 p^3+15 p^4+9 p^5+3 p^6 \cr \noalign{\hrule}
 3 & 4 & 4 & 12+10 p+18 p^2+14 p^3+18 p^4+12 p^5+12 p^6+4 p^7 \cr \noalign{\hrule}
 4 & 4 & 4 & 13+11 p+20 p^2+16 p^3+21 p^4+15 p^5+16 p^6+8 p^7+5 p^8 \cr} \hrule}
\]

\vspace{2cm}

\noindent\textbf{Mario Hampejs }\\
NuHAG, Faculty of Mathematics, Universit\"at Wien\\
Nordbergstra{\ss}e 15, A-1090 Wien \\ Austria\\
{\tt mario.hampejs@univie.ac.at}\\

\noindent\textbf{L\'aszl\'o T\'oth}\\
Institute of Mathematics, Universit\"at f\"ur Bodenkultur \\
Gregor Mendel-Stra{\ss}e 33, A-1180 Wien \\ Austria \\ and \\
Department of Mathematics, University of P\'ecs \\ Ifj\'us\'ag u. 6,
H-7624 P\'ecs \\ Hungary \\ {\tt ltoth@gamma.ttk.pte.hu}


\vspace{-3ex}
\begin{thebibliography}{99}
\setlength{\parskip}{-3pt}\vspace{-2ex}

\bibitem{Aig2007} \textbf{Aigner, M.,} \textit{A Course in Enumeration}, Graduate
Texts in Mathematics 238, Springer, 2007.

\bibitem{Bho1996} \textbf{Bhowmik, G.,} Evaluation of divisor functions of matrices,
\textit{Acta Arith.}, \textbf{74} (1996), 155--159.

\bibitem{Bir1934} \textbf{Birkhoff, G.} Subgroups of Abelian groups, \textit{Proc.
London Math. Soc.}, (2) {\bf 38} (1934--35), 385--401.

\bibitem{But1994} \textbf{Butler, M.~L.,} \textit{Subgroup Lattices and Symmetric Functions},
Mem. Amer. Math. Soc. \textbf{112} (1994), no. 539.

\bibitem{Cal2004} \textbf{C\u{a}lug\u{a}reanu, G.,} The total number of subgroups
of a finite abelian group, \textit{Sci. Math. Jpn.}, \textbf{60}
(2004), 157--167.

\bibitem{Del1948} \textbf{Delsarte, S.,} Fonctions de M\"obius sur les groupes abeliens
finis, \textit{Annals of Math.}, \textbf{49} (1948), 600--609.

\bibitem{Dyu1948} \textbf{Dyubyuk, P.~E.,} On the number of subgroups of a finite
abelian group, \textit{Izv. Akad. Nauk SSSR Ser. Mat.}, \textbf{12}
(1948), 351--378.

\bibitem{HHTW2012} \textbf{Hampejs, M., Holighaus, N., T\'oth, L. and Wiesmeyr, C.,} On the subgroups
of the group $\Z_m \times \Z_n$, submitted, \\ {\tt
http://arxiv.org/abs/1211.1797}

\bibitem{Hux2003} \textbf{Huxley, M. N.,} Exponential sums and
lattice points III., \textit{Proc. London Math. Soc.}, {\bf 87} (2003),
591--609.

\bibitem{Pet2011} \textbf{Petrillo, J.,} Counting subgroups in a direct product of finite cyclic grups,
\textit{College Math J.}, {\bf 42} (2011), 215--222.

\bibitem{Rem1931} \textbf{Remak, R.,} \"Uber Untergruppen direkter Produkte von drei Faktoren,
\textit{J. Reine Angew. Math.}, \textbf{166} (1931), 65--100.

\bibitem{Sch1994} \textbf{Schmidt, R.,} \textit{Subgroup Lattices of Groups}, de Gruyter Expositions
in Mathematics 14, de Gruyter, Berlin, 1994.

\bibitem{OEIS} \textbf{Sloane, N.~J.~A.,} \textit{On-Line Encyclopedia of Integer Sequences (OEIS)},
\\ http://oeis.org.

\bibitem{Suz1951} \textbf{Suzuki, M.,} On the lattice of subgroups of finite
groups, \textit{Trans. Amer. Math. Soc.}, \textbf{\bf 70} (1951),
345--371.

\bibitem{Tar2007} \textbf{T\u{a}rn\u{a}uceanu, M.,} A new method of proving some classical
theorems of abelian groups, \textit{Southeast Asian Bull. Math.},
\textbf{31} (2007), 1191--1203.

\bibitem{Tar2010} \textbf{T\u{a}rn\u{a}uceanu, M.,} An arithmetic method of counting the
subgroups of a finite abelian group, \textit{Bull. Math. Soc. Sci.
Math. Roumanie (N.S.)}, \textbf{53(101)} (2010), 373--386.

\bibitem{Tot2010} \textbf{T\'oth, L.,} A survey of gcd-sum functions, \textit{J. Integer
Sequences}, \textbf{13} (2010), Article 10.8.1, 23 pp.

\bibitem{Tot2011} \textbf{T\'oth, L.,} Menon's identity and arithmetical sums
representing functions of several variables, \textit{Rend. Sem. Mat.
Univ. Politec. Torino}, \textbf{69} (2011), 97--110.

\bibitem{Tot2012} \textbf{T\'oth, L.,} On the number of cyclic subgroups of a finite Abelian
group, \textit{Bull. Math. Soc. Sci. Math. Roumanie (N.S.)},
\textbf{55(103)} (2012), 423--428.

\bibitem{Yeh1948} \textbf{Yeh, Y.,} On prime power abelian groups, {\it Bull. Amer.
Math. Soc.}, \textbf{54} (1948), 323--327.
\end{thebibliography}
\end{document}